\documentclass[leqno,a4paper,11pt]{amsart}

\RequirePackage[dvips]{hyperref}
\usepackage{amsmath,amsthm,amssymb}

\numberwithin{equation}{section}
\theoremstyle{plain}
\newtheorem{theo}{Theorem}[section]
\newtheorem{lemm}{Lemma}[section]

\theoremstyle{remark}
\newtheorem{rema}{Remark}[section]

\usepackage[latin1]{inputenc}
\usepackage{aeguill}
\usepackage{enumerate,amssymb}

\long\def\xcom#1{}

\newcommand{\ladate}{\space\the\day\ \ifcase\month\or janvier\or f\'evrier
\or mars\or avril\or mai\or juin\or juillet\or aout\or septembre
\or octobre\or novembre\or d\'ecembre\fi
\ {\oldstyle\the\year}}

\def\alphabet#1{\ifcase#1\or a\or b\or c\or d\or e\or f\or g\or 
h\or i\or j\or k\or l\or m\or n\or o\or p\or q\or r\or s\or t\or 
u\or v\or w\or x\or y\or z\fi}

\newcommand{\R}{\mathbb{R}}

\newcommand{\unsur}[1]{{\frac{1}{#1}}}
\def\un#1{{{\,\mathbf{1}}_{({#1})}}}
\def\valabs#1{{\left\vert {#1} \right \vert}}

\def\etp#1{{\left ( {#1} \right )}}
\def\etc#1{{\left [ {#1} \right ]}}
\def\crochet#1{{\left < {#1} \right >}}

\def\norme#1{{\left \Vert #1 \right \Vert}}
\def\ens#1{{\left\{#1\right\}}}

\def\dessus#1#2{\mathord{\mathop{\kern 0pt #2}\limits^#1}}

\newcommand{\bit}{\begin{itemize}}
\newcommand{\eit}{\end{itemize}}
\newcommand{\ben}{\begin{enumerate}}
\newcommand{\een}{\end{enumerate}}

\newcounter{moncompteur}

\newenvironment{myenumerate}%
{\begin{list}{\arabic{moncompteur}. }{\usecounter{moncompteur}%
\setlength{\leftmargin}{0pt}%
\setlength{\labelwidth}{0pt}%
\setlength{\listparindent}{0pt}%
\setlength{\labelsep}{0pt}}}%
{\end{list}}

\def\bmen{\begin{myenumerate}}
\def\emen{\end{myenumerate}}

\newcommand{\intot}{{\int_0^t}}

\newcommand{\undemi}{\frac{1}{2}}

\def\l{\lambda}

\newcommand{\Arond}{{\mathcal A}}

\newcommand{\Hrond}{{\mathcal H}}

\newcommand{\Lrond}{{\mathcal L}}
\newcommand{\Mrond}{{\mathcal M}}

\newcommand{\Srond}{{\mathcal S}}

\newcommand{\Urond}{{\mathcal U}}

\usepackage{fancybox}

{\par\begin{Sbox}\begin{minipage}{\textwidth}}%
{\end{minipage}\end{Sbox}\par\smallskip\shadowbox{\TheSbox}}

\newcommand{\Lie}{\mathfrak{L}}

\newcommand{\ppun}{\partial_{p_1}}

\newcommand{\pqun}{\partial_{q_1}}

\newcommand{\prun}{\partial_{r_1}}
\DeclareMathOperator{\signe}{sign}
\begin{document}

\title[Chain of Oscillators]{Existence and uniqueness of an invariant measure for a chain of oscillators in contact with two heat baths}

\author{{Philippe} {Carmona}}


\address{Philippe Carmona\\
Laboratoire Jean Leray, UMR 6629,
Universit{\'e} de Nantes, BP 92208\\
F-44322 Nantes Cedex 03\\
E-mail : philippe.carmona@math.univ-nantes.fr\\
\url{http://www.math.sciences.univ-nantes.fr/~carmona}
}

\begin{abstract}
In this note we consider a chain of $N$ oscillators, whose ends are in contact with two heat baths at different temperatures. Our main result is the exponential convergence to the unique invariant probability measure (the stationary state).
We use the Lyapunov's function technique of Rey-Bellet and coauthors~\cite{MR1991548,MR1915300,MR1889227,MR1799873,MR1685893,MR1705589}, with different model of heat baths, and adapt these techniques to two new case recently considered in the literature by Bernardin and Olla~\cite{MR2185330} and Lefevere and Schenkel~\cite{1742-5468-2006-02-L02001}.
\end{abstract}

\subjclass{Primary 60K35, 60J65, 82C99}
\keywords{
{Hamiltonian chain} ;
{harmonic oscillator} ;
{heat bath} ;
{invariant measure}
}

\date{\today}\vfuzz=4pt \hfuzz=40pt
\maketitle

\section{Introduction}

We consider a chain of coupled oscillators whose  dynamic is described by the
Hamiltonian 

\begin{equation*}
  H = \sum_{1\le i\le N} \undemi p_i^2 + V(q_i) + \sum_{1\le i\le N-1} U(q_{i+1} -q_i)\,,
\end{equation*}


where $q_i$ and $p_i$ are the position and the momentum of the oscillators. We assume that the potentials $U,V$ have the following properties:

\textbf{(H1) : Growth at infinity} $U$ and $V$ are $C^\infty$ and there exists real constants $l\ge 1$, $k\ge 2$, $k\ge l$, $a_k>0,b_l >0$ such that:
\begin{align*}
  \lim_{\l \to +\infty} \l^{-k} U(\l x) = a_k \valabs{x}^k\,,& \lim_{\l \to +\infty} \l^{1-k} U'(\l x) = ka_k \valabs{x}^{k-1} \signe(x)\\
 \lim_{\l \to +\infty} \l^{-l} V(\l x) = b_l \valabs{x}^l\,,& \lim_{\l \to +\infty} \l^{1-l} V'(\l x) = l a_l \valabs{x}^{l-1} \signe(x)\\
\end{align*}

\textbf{(H2) : Non degeneracy} The interaction potential is non degenerate : for any $q \in \R$ there exists $m=m(q)\ge 2$ such that $\partial^m U(q) \neq 0$.

\smallskip

We want to stress the fact that $k\ge l$, that is \emph{ near infinity the interaction potential $U$ dominates the pinning potential $V$}. Typical examples of such potentials are:
\begin{itemize}
\item The harmonic chain : $U(x) = \undemi x^2$, $V = \alpha x^2$, with $\alpha \ge 0$.
\item The Fermi-Pasta-Ulam chain : $U(x) = \undemi x^2 + \frac{x^4}{4}$, $V(x) = \alpha x^4$, with $\alpha \ge 0$.
\item The harmonic coupling with weak pinning : $U(x) = \undemi x^2$, $V(x)=(1+x^2)^{\undemi}$.
\end{itemize}

The two ends of the chain are in contact with heat baths at temperatures $T_1=T_l \ge T_r = T_N$. The interaction between the heat bath and the particle is modeled by a Langevin process (an Ornstein Uhlenbeck process) at the corresponding temperature. Therefore the dynamics are described by the following system of stochastic differential equations:

\begin{equation}\left\{
 \begin{aligned}
   dq_i(t) &= p_i(t)\, dt\\
   dp_i(t) &= (-\partial_{q_i}H -\undemi p_i)\, dt + \sqrt{T_i} dB_i(t) & (i=1,N)\\
dp_i(t) &= (-\partial_{q_i}H)\, dt & (2\le i\le N-1)\,.
\end{aligned}  \right.
\tag{\text{$S$}}
\end{equation}


where $B_1,B_N$ denote two independent standard Brownian motions.

In the \emph{unpinned case}, that is $V=0$, the dynamics are translation invariant (ie under the action $q_i \to q_i + C$), and there is no hope of finding an invariant probability measure. It is thus natural to consider the \emph{interdistances}
$r_i = q_{i+1} -q_i$, for $1 \le i\le N-1$, and the system

\begin{equation}\left\{
 \begin{aligned}
   dr_i(t) &= (p_{i+1} -p_i)\, dt\\
   dp_1(t) &= (U'(r_1) -\undemi p_1)\, dt + \sqrt{T_1} dB_1(t) & \\
   dp_N(t) &= (-U'(r_{N-1}) -\undemi p_N)\, dt + \sqrt{T_N} dB_N(t) & \\
   dp_i(t) &= (U'(r_i) -U'(r_{i-1}))\, dt & \text{\small $(2\le i\le N-1)$}
\end{aligned}  \right.
\tag{\text{$S'$}}
\end{equation}

Observe that the state space of the stochastic system is $\chi = \R^{2N}$ in the pinned case, and $\chi=\R^{2N-1}$ in the unpinned case.

\medskip
To state the main result, we need to introduce the Lyapunov function $W = e^{\theta H}$, for a $\theta >0$ to be chosen later, and the corresponding weighted Banach space
$$ \Hrond = \ens{ f:\chi \to \R \text{ continuous, }\norme{f}_W < +\infty}\,,\qquad \text{with}\quad \norme{f}_W = \sup_x \frac{\valabs{f(x)}}{W(x)}$$


\begin{theo}\label{sec:introduction:thm}
There exists $\theta>0$ such that:
\begin{enumerate}[(i)]
  \item For every starting point, the system $(S)$ (resp. $(S')$) has a unique solution, defined for all times $t\ge 0$.
\item The corresponding semigroup $(T_t)_{t\ge 0}$ has a smooth density.
\item The Markov system has a unique invariant probability measure $\pi$, which is absolutely continuous with a smooth density.
\item The semigroup converges exponentially fast to the invariant measure (and is therefore ergodic). More precisely, there exists a constant $C>0$ such that
$$\valabs{T_tf(x) -\pi(f)} \le C \norme{f}_W\,e^{-C t} W(x) \qquad(\forall x\in\chi, f \in \Hrond)$$
\end{enumerate}
  
\end{theo}

The proofs we give are similar to the original ones, as they are stated in the survey paper~\cite{reybelletopen}, or in  the original papers by Rey-Bellet and Thomas~\cite{MR1915300,MR1889227}, and Eckmann, Pillet and Rey-Bellet~\cite{MR1685893,MR1705589}. However we believe that this note is an interesting contribution to the subject, because of the following:

\begin{itemize}
\item The model of heat baths is the simplest possible, an Ornstein-Uhlenbeck process. We do not introduce the auxiliary variables of \cite{MR1915300,MR1889227,MR1685893,MR1705589,MR1764365}. However, we do not provide a nice physical interpretation of such a model of heat bath.
\item The assumptions on the potentials are slightly relaxed : we do not require that the pinning potential to be at least quadratic, that is $l\ge 2$. Observe that we do not reach the level of generality of Eckmann and Hairer \cite{MR1764365}, where only growth bounds on $U,V$ and their derivatives are required.
\item In  \cite{MR1915300,MR1889227,MR1685893,MR1705589,MR1764365}, the authors use a scaling argument based on two ingredients :
  \begin{itemize}
  \item Continuity of solutions of ordinary differential equations with respect to initial conditions and coefficients.
  \item A tracking lemma which says that at high energy the system is nearly deterministic (one can almost forget the heat bath influence).
  \end{itemize}
We replaced both arguments by a single one : continuity of stochastic differential equations with respect to initial conditions and coefficients. This is not only a technical shortcut, but also enables the study of models of chain where in addition to heat baths, you have a specific noise that you need to keep even at high energy since it is a feature of the model.
\item The first of these models was introduced by by Bernardin and Olla~\cite{MR2185330}. There is, in addition to heat baths,  a random exchange of momentum between neighbouring atoms. The model at infinite energy is not deterministic but stochastic (see the system~\ref{bolla:sysinfini} on page \pageref{bolla:sysinfini}).
\item The second model we study is 
the chain introduced by Lefevere and Schenkel~\cite{1742-5468-2006-02-L02001}, where the noise is highly non standard : in Fourier coordinates the momenta are in contact with heat baths at the same temperature $T$ and positions are also coupled to the heat baths. It is not a straightforward application of results of Rey-Bellet and coauthors since we need to show that the energy dissipation on the momenta is the leading term in the upper bound of \eqref{eq:lefevere:majsg}. Alas we are not able to generalize the results to non quadratic pinning potentials $V$. In particular we do not know how assumptions of the type $C_1 \le V'' \le C_2$ could be translated in Fourier coordinates.
\end{itemize}


\section{Non explosion}
A priori, if we consider the stochastic differential equation
\begin{equation}
  \label{nonexp:eq:1}
  dX_t = b(X_t)\, dt + \sigma(X_t)\, dB_t
\tag{\text{$e_{\sigma,b}$}}
\end{equation}

with locally Lipschitz coefficients $b:\R^n \to \R^n$, $\sigma : \R^n \to \Mrond_{n\times m}$, $B$ a Brownian motion in $\R^m$, we can only be sure, given a starting point $x$, of the existence of a local solution, defined up to an explosion random time $\zeta(\omega)$. However, if there exists a Lyapunov function $W$, ie a function such that $W\ge 1$ and the level sets $\ens{W \le A}$ are compacts, that satisfies $L W \le C W$ for a constant $C>0$, with $L$ the formal generator of the diffusion $X$, that is the second order differential operator 
$$ L = \sum_i b_i(x) \partial_{x_i} + \undemi \sum_{i,j} a_{ij}(x) \partial^2_{xi\,x_j} \qquad (a=\sigma \sigma^t)\,,$$
then (see Theorem 5.9 of~\cite{reybellet03}) the lifetime is $\zeta=+\infty$, that is $X$ is defined for all times $t\ge 0$, with semigroup $(T_t)_{t\ge 0}$ such that :
$$ T_t W(x) = P^x\etp{W(X_t)} \le W(x)$$
This ensures that $(T_t)_{t\ge 0}$ can be extended to a strongly continuous semigroup on $\Hrond$.

For the pinned chain of oscillators, we have 
$$ L = \Arond + \Lrond_R\,,\quad\text{with}\quad \Arond = \sum_i \partial_{p_i}H \partial_{q_i} - \partial_{q_i}H \partial_{p_i}$$
the Hamiltonian generator and $\Lrond_R$ the heat baths (reservoirs) generator
$$ \Lrond_R = \undemi \sum_{i=1,N} (-p_i \partial_{p_i} + T_i \partial^2_{p_ip_i})\,.$$

The energy $H$ is conserved by the Hamiltonian dynamics : $\Arond H =0$. Since $\Arond$ is a first order differential operator 
$$ \Arond W = e^{\theta H} \theta \Arond(H) = 0\,.$$
We shall make the {\emph{assumption}}
$$ 0 < \theta < \unsur{max(T_1,T_N)}$$
Then, 
\begin{gather*}
  \Lrond_R H = \undemi ( T_1 - p_1^2 + T_N -p_N^2) \\
L W = \Lrond_R W = \undemi \theta W \etp{(T_1 + T_N) - (p_1^2 + p_N^2) + \theta(T_1 p_1^2 + T_N p_N^2)} \le \undemi \theta (T_1 + T_N) \, W
\end{gather*}

For the unpinned chain, we obtain exactly the same upper bound, since we still have $\Arond H = \Arond W = 0$.

Hence, part (i) of Theorem~\ref{sec:introduction:thm} is established.

Observe that since $V(x)$ and $U(x)$ go to infinity as $\valabs{x}\to +\infty$,  $H$ is  bounded below, therefore we need to multiply $e^{\theta H}$ by a constant $C$ in order to have $W=Ce^{\theta H}\ge 1$. We shall forget about this constant in the sequel since it will not change any proof.

\section{Hypoellipticity}
The generator of the pinned system can be written

$$ L = X_0 + X_1^2 + X_N^2\,,$$
where $X_i$ are first order differential operators
$$ X_0 = \Arond - \undemi(p_1 \partial_{p_1} + p_N \partial_{p_N})\,,\quad X_i = \sqrt{\frac{T_i}{2}} \partial_{p_i} \quad(i=1,N)\,,$$
and $\Arond = \sum_i \partial_{p_i}H \partial_{q_i} - \partial_{q_i}H \partial_{p_i}$ is the Hamiltonian generator. 

Let $\Lie$ be the Lie algebra generated by the vector fields
$$ \ens{X_i}_{i\ge 1}, \ens{\etc{X_i,X_j}}_{0\le i,j}, \ens{\etc{X_i,\etc{X_j,X_k}}}_{0\le i,j,k}\,, \ldots$$
Assume H\"ormander's hypoellipticity condition, that is that $\Lie$ has full rank  at every point $x$. Then (see e.g.  Corollary~7.2  of~\cite{reybellet03} and the references therein, or H\"ormander's original paper~\cite{MR0222474}), 
\begin{itemize}
\item The semigroup $T_t$ has a smooth density, i.e.
$$ T_t f(x) = \int p_t(x,y) f(y)\, dy$$
with $(t,x,y) \to p_t(x,y)$ smooth.
\item The semigroup $T_t$ is Strong Feller : for $t>0$, $T_t$  sends measurable bounded functions into continuous functions.
\item The invariant measures, if they exist, also have a smooth density.
\end{itemize}

Hence, we only need to verify that $\Lie$ has full rank, in both pinned and unpinned cases.

\textbf{Pinned case} $\Lie$ contains $X_1,X_N$ so it contains $\partial_{p_1}$ and $\partial_{p_n}$. But,
$$ \etc{\ppun,X_0} = -\undemi \ppun + \etc{\ppun,\Arond} = -\undemi \ppun + \pqun $$
Therefore, $\pqun \in \Lie$. Since,
$$ \etc{\pqun,X_0}=\etc{\pqun,\Arond} = - \sum_i \partial_{q_1q_i H}\, \partial_{p_i} = -\partial_{q_1q_1}H \ppun -\partial_{q_1q_2} H \, \partial_{p_2}$$
we have $\partial_{q_1q_2} H \, \partial_{p_2} \in \Lie$. We iterate this procedure, by computing $\etc{\pqun,\partial_{q_1q_2} H \, \partial_{p_2}}$ until we obtain that $\partial_{q_1^{m-1} q_2} H \partial_{p_2} =(-1)^{m-1} \partial^m U(q_2-q_1)\, \partial_{p_2} \in \Lie$. By the non degeneracy assumption, this implies  $\partial_{p_2} \in \Lie$ and $\etc{\partial_{p_2},X_0}=q_2 \in \Lie$. By induction, we establish that for all i, $\partial_{p_i}$ and $\partial_{q_i}$ are in $\Lie$.

\medskip

\textbf{Unpinned case}
Only the Hamiltonian generator is different :
$$\Arond = \sum_{1\le i\le N-1} (p_{i+1}-p_i) \partial_{r_i} + \sum_{2\le i\le N-1} (U'(r_i) -U'(r_{i-1})) \partial_{p_i} + U'(r_1) p_1 - U'(r_{N-1}) \partial_{p_N}\,.$$
By assumption $\ppun \in \Lie$ and
$$ \etc{\ppun, X_0}= -\undemi \ppun - \prun$$
Hence $\prun \in \Lie$ and
$$ \etc{\prun , X_0} = -U''(r_1) \partial_{p_2} + U''(r_1) \ppun$$
Consequently, $U''(r_1) \partial_{p_2}\in\Lie$ and
$$ \etc{\prun,U''(r_1) \partial_{p_2}} = \partial^3_{r_1^3} U(r_1) \partial_{p_2}\in\Lie$$
and we obtain by induction that $\partial^m_{r_1^m} U(r_1) \partial_{p_2}\in\Lie$. Therefore, $\partial_{p_2}\in\Lie$ and since
$$ \etc{\partial_{p_2},X_0} = \prun \partial_{r_2}$$
we have $\partial_{r_2} \in \Lie$ and by induction we show that $\Lie$ has full rank $2N-1$.

\begin{rema}
  Observe that we only need one heat bath to ensure hypoellipticity.
\end{rema}

\section{The Control Problem}
The way chosen to establish uniqueness of the invariant measure, if it exists, is to prove that the semigroup $(T_t)_{t\ge 0}$ is \emph{irreducible}, that is for every $x$, every $t>0$, every non empty open set $A$: $T_t(x,A)>0$. The most frequently used tool to establish irreducibility is Stroock-Varadhan's support theorem\cite{MR0400425}: let $X$ be the solution of the stochastic differential equation, in the Stratonovitch sense,
$$ \partial X_t = b(X_t) \partial_t + \sigma(X_t) \partial B_t$$
Let $\Urond=\ens{u: \R_+ \to \R^m: u(0)=0, u'\in \Lrond^2(0,t) \forall t}$ be the set of controls lying in the Cameron-Martin space. Let $(C_{b,\sigma,u})$ be the controlled ordinary differential equation
\begin{equation}
  \label{control:eq:1}
  \frac{dx}{dt} = b(x) + \sigma(x) u(x)
\end{equation}
Let $x(x_0,u,t)$ be the maximal integral curve of the controlled equation, passing through $x_0$ at time $t=0$. Let $A(x_0,t)$ (resp $A^*(x_0,t)$ ) be the accessibility set (resp the strong accessibility set) at time $t>0$ starting from $x_0$:
$$ A(x,t) = \ens{x(x_0,u,s) : u\in\Urond, 0\le s\le t}\,,\quad
 A^*(x,t) = \ens{x(x_0,u,t) : u\in\Urond}$$
The the support theorem states that the support of the measure $T_t(x_0,.)$ is $cl(A^*(x_0,t))$ the closure of the strong accessibility set.

In the general case, establishing that $cl(A^*(x_0,t))=\chi(=\R^{2N} \text{ or } \R^{2N-1})$ is not a trivial thing, and we refer to Eckmann,Pillet and Rey-Bellet\cite{MR1705589}, Theorem~3.2 for a proof.

However, for the harmonic chain, $U(x)=\undemi x^2$ and $V(x)=\undemi \alpha x^2$,  we can provide a simple proof, that we learnt from François Laudenbach. Indeed, the control system is linear. First since the matrix $\sigma(x)$ is constant, there is no difference between Itô and Stratonovitch integrals. Second, let us write the control system, for the pinned case: we incorporate the constants $\sqrt{T_i}$ into the two controls $u_1,u_N$ so that we obtain

\begin{equation}\left\{
 \begin{aligned}
   dq_i(t) &= p_i(t)\, dt\\
   dp_i(t) &= (-\partial_{q_i}H -\undemi p_i)\, dt + u_i & (i=1,N)\\
dp_i(t) &= (-\partial_{q_i}H)\, dt & (2\le i\le N-1)\,.
\end{aligned}  \right.
\tag{\text{$C$}}
\end{equation}

This can be written, if $x=(p,q) \in\R^{2N}$ as
$$ \frac{dx}{dt} = A x + B u(t)$$
with $A$ a $n\times n$ matrix and $B$ a $n\times 2$ matrix ($n=2N$). Kalman's criterion states that for every $t>0$ and $x\in\chi=\R^n$, we have $A^*(x_0,t)=\R^n$ as soon as the smallest vector space $\Srond$ containing the image of $B$, and stable by $A$, is the whole space $\chi=\R^n$ (see e.g. Wonham's\cite{MR569358}).

Here the image of $B$ is spanned by $\ens{e_{N+1},e_{2N}}$ where the $e_i$ are the canonical base of $\R^n$. Assume for example that $V(q)=\undemi \alpha q^2$ and $U(q)=\undemi q^2$. then the control system is

\begin{equation}\left\{
 \begin{aligned}
   \frac{dq_i}{dt} &= p_i(t)\\
\frac{dp_1}{dt} &= q_2 - q_1 + \alpha q_1+ u_1 \\
\frac{dp_N}{dt} &= -(q_N - q_{N-1}) + \alpha q_N + u_2 \\
\frac{dp_i}{dt} &=  (q_{i+1}+q_{i-1} - 2 q_i) + \alpha q_i& (2\le i\le N-1)\,.
\end{aligned}  \right.
\tag{\text{$C$}}
\end{equation}

The matrix $A$ can be written as blocks
$$ A =\begin{pmatrix} 0 & \tilde{A} \\ I & 0 \end{pmatrix}$$
Therefore we only need to check that the smallest vector space containing $\tilde{e}_1$ and $\tilde{e_N}$, and stable by $\tilde{A}$ is $\R^N$ itself. This is obvious since we have 

$$ \tilde{A}=
\begin{pmatrix}
  -1+\alpha & 1 & 0 & 0 &\ldots & 0 \\
1 & -2 + \alpha & 1 & 0 & \ldots & 0 \\
\ldots & \ldots & \ldots & \ldots & \ldots & \ldots \\
0 &0 & 0 & 0 & \ldots & -1 + \alpha 
\end{pmatrix}
$$

\section{Scaling}
\label{sec:scaling}

\newcommand{\ale}{E^{\unsur{k} -\undemi}}
\newcommand{\alen}{E_n^{\unsur{k} -\undemi}}

We shall scale both positions $q_i$ and moments $p_i$ and use continuity with respect to parameters, for solution of stochastic differential equations, to prove that the system is attracted to the compact sets where the energy $H$ stays bounded. More precisely, given $X(t)=(q(t),p(t))$ solution of $(S)$ starting from $x=(q(0), p(0))$, and $E>0$, we shall consider

$$ X^E(t) = \left\{
 \begin{aligned}
   p_i^E(t) &=E^{-\undemi} p_i(\ale t) \\
  q_i^E(t) &= E^{-\unsur{k}} q_i(\ale t)
\end{aligned}  \right. $$
Of course $\frac{dq^E_i}{dt} = p^E_i$ and for $2\le i\le N-1$, we have 
$$ \frac{dp_i^E}{dt} = - \partial_{q_i} H_E(q^E(t)) $$
with the Hamiltonian
$$ H_E(p,q) = \sum_i \undemi p_i^2 + \unsur{E} V(E^{\unsur{k}} q_i) + \sum_{1\le i\le N-1} \unsur{E} U(E^{\unsur{k}} (q_{i+1} -q_i))\,.$$
For $i=1,N$ we obtain
$$ dp_i^E(t) = E^{-\undemi} p(0) - \intot (\ale \undemi p^E_i(s) + \partial_{q_i} H (q^E(s)))\, ds + E^{\unsur{2k} - \frac{3}{4}} \sqrt{T_i} B^E_i(t)$$
where $B^E(t) = E^{\unsur{4} - \unsur{2k}} B(\ale t)$ is a standard Brownian motion. In a nutshell $X^E$ is the solution, starting from $x^E=(E^{-\undemi} p(0),E^{-\unsur{k}} q(0))$ of the stochastic differential system

\begin{equation}\left\{
 \begin{aligned}
   dq_i(t) &= p_i(t)\, dt& \\
   dp_i(t) &= (-\partial_{q_i}H_E -\frac{\ale}{2} p_i)\, dt +  E^{\unsur{2k} - \frac{3}{4}} \sqrt{T_i} dB_i(t) & (i=1,N)\\
dp_i(t) &= (-\partial_{q_i}H_E)\, dt & \text{\small $(2\le i\le N-1)$}\,.\\
\end{aligned}  \right.
\tag{\text{$S^E$}}
\end{equation}
Observe that the scaling is such that if $H(x)=E$ then $H_E(x^E)=1$.

The main ingredient of the proof is the convergence of a solution of $(S^E)$ to a solution of the limit system, where the noise has disappeared,

\begin{equation}\left\{
 \begin{aligned}
   dq_i(t) &= p_i(t)\, dt& \\
   dp_i(t) &= (-\partial_{q_i}H_\infty -\undemi \un{k=2} p_i)\, dt  & (i=1,N)\\
dp_i(t) &= (-\partial_{q_i}H_\infty)\, dt & \text{\small $(2\le i\le N-1)$}\,.\\
\end{aligned}  \right.
\tag{\text{$S^\infty$}}
\end{equation}
with $$ H_\infty(p,q) = \sum_i \undemi p_i^2 + \un{k=l} b_k \valabs{q_i}^k + \sum_{1\le  i\le N-1} a_k \valabs{q_{i+1}-q_i}^k.$$

\medskip

We shall establish now that the system $(S^\infty)$ is non degenerate, that is that the chain remains still when both ends are still.

\begin{lemm}\label{sec:scaling:lem}
  Assume $l=k$. Let $(q(t),p(t))$ be a solution of $(S^\infty)$, starting from $x$ such that $H_\infty(x)=1$. Then, for any $\tau >0$,
$$ \int_0^\tau (p_1^2(s) + p_N ^2(s))\, ds >0\,.$$
\end{lemm}
\begin{proof}
\xcom{
  If $degree(V)=degree(U)$, that is $k=l$, the $\ens{x : H_\infty(x) =1}$ is compact, the second part follows immediately from the first and from the continuity of solution of ordinary differential equations with respect to the starting point.
}

Assume that $\int_0^\tau (p_1^2(s) + p_N ^2(s))\, ds =0$ ; we are going to obtain a contradiction. Let $I=\etc{0,\tau}$. Then, by assumption, $p_1 = 0 $ on $I$ so $\frac{dq_1}{dt}= p_1 =0$ and $q_1 = c_1$ is constant on $I$. We have, on $I$
$$ 0 = \frac{dp_1}{dt} = - (\un{k=2} \undemi p_1 + \partial_{q_1} H_\infty)$$
therefore
$$0 = b_k \valabs{c_1}^{k-1} \signe(c_1) + a_k \valabs{q_2 -c_1}^{k-1}\signe(q_2-c_1)$$
hence $q_2=c_2$ is constant on $I$, and by induction we obtain that $p_i=0$ and $q_i=c_i$ are constant on $I$. If we write $H_\infty(p,q)=\sum_{i} \undemi p_i^2 + g(q)$ we obtain that $c=(c_1, \ldots, c_N)$ is a solution of $\nabla g(c)=0$. Since $g$ is convex ($k$ is even), $c$ is a global minimum of $g$, which is absurd because $g(0)=0$ and $g(c)=H_\infty(0,c)=1$ by assumption.
\end{proof}

\begin{rema}

  For the unpinned case we obtain exactly the same result, for the limit system 
\begin{equation}\left\{
 \begin{aligned}
   dr_i(t) &= (p_{i+1} -p_i)\, dt& \\
   dp_1(t) &= (k a_k \valabs{r_1}^{k-1}-\undemi \un{k=2} p_1)\, dt  & \\
   dp_N(t) &= (-k a_k \valabs{r_{N-1}}^{k-1}-\undemi \un{k=2} p_N)\, dt  & \\
dp_i(t) &= k a_k (\valabs{r_i}^{k-1} -\valabs{r_{i-1}}^{k-1})\, dt & \text{\small $(2\le i\le N-1)$}\,.\\
\end{aligned}  \right.
\tag{\text{${S'}^\infty$}}
\end{equation}
If $H_\infty(x)=1$ and $X$  a solution of $({S'}^\infty)$, starting from $x$, then $\int_0^\tau (p_1^2(s) + p_N ^2(s))\, ds >0$.
\end{rema}

From the preceding Lemma, we shall derive the asymptotic result:


\begin{lemm}\label{sec:scaling-1}
  Assume that $(X_n(t))_{t\ge 0}$ is a solution of $(S)$ starting from $x_n$, with $E_n=H(x_n)\to +\infty$. Then there exists a subsequence, $(x_{n_k})_k$ such that for any $C>0, t_0>0$, 
$$ \lim_{k \to +\infty} P^{x_{n_k}}\etc{\exp\etp{- C \int_0^{t_0} (p_1^2(s) + p_N ^2(s))\, ds}} = 0\,.$$
\end{lemm}


\begin{proof}
  The process $(X^{E_n}_n (t))_{t\ge 0}$ is a solution of $(S^{E_n})$ starting from $x^{(n)}$ with $H_{E_n}(x^{(n)}) = \unsur{E_n} H(x_n)= 1$. We build, on the same filtered probability space, that is with the same Brownian $B$, a family of processes $({X}^{(n)}(t))_{t\ge 0}$, solution of $(S^{E_n})$ starting from $x^{(n)}$.

Assume first that $k=l$. Since $H_{E_n}(x^{(n)})=1$, the sequence $x_n$ remains in a compact set and we can extract a converging subsequence, which we shall still denote by $x^{(n)}$. Hence we have $x^{(n)}\to x^{(\infty)}$ and by continuity $H_\infty( x^{(\infty)})=1$. 

Thanks to the continuity of solutions of stochastic differential equations with respect to both parameters and starting points (see e.g. Bahlali, Mezerdi and Ouknine~\cite{MR1655150}) we have, for any $\tau>0$:
$$ \lim_{n \to +\infty} P\etc{\sup_{t\le \tau} (X^{(n)}(t) - X^{(\infty)}(t))^2} =0$$
where $X^{(\infty)}$ is the solution of $(S^\infty)$ starting from $x^{(\infty)}$.

Therefore, for any $A>0$, we infer from Lemma~\ref{sec:scaling:lem} that if
$c_\tau =  \int_0^\tau  (p_1^2(s) + p_N ^2(s))\, ds$ for $X^{(\infty)}$, then
\begin{equation*}
  \lim_{n\to +\infty} P^{x^{(n)}}\etc{e^{-A \int_0^\tau  (p_1^{E_n}(s)^2 + p_N^{E_n}(s)^2) \, ds}} 
= P^{x^{(\infty)}}\etc{e^{-A \int_0^\tau  (p_1^2(s) + p_N ^2(s))\, ds}} 
=e^{-A c_\tau}
\end{equation*}

Observe now that
$$ \int_0^{t_0} (p_1^2(s) + p_N^2(s))\, ds = E_n^{\unsur{k} + \undemi}\int_0^{t_0 E_n^{\undemi -\unsur{k}}}  (p_1^{E_n}(s)^2 + p_N^{E_n}(s)^2)\, ds$$
Since $k\ge 2$, choosing $\tau=t_0$, we see that  there exists $n_0$ such that $t_0 E_n^{\undemi -\unsur{k}} \ge \tau$ for $n\ge n_0$. Finally,

\begin{align*}
\limsup_{n\to +\infty} P^{x_n}\etc{e^{-C  \int_0^{t_0} (p_1^2(s) + p_N ^2(s))\, ds}} & \le \limsup_{n\to +\infty} P^{x^{(n)}}\etc{e^{-C  E_n^{\unsur{k} + \undemi}\int_0^{\tau} (p_1^{E_n}(s)^2 + p_N^{E_n}(s)^2)\, ds}}
\\ & = 0\,.
\end{align*}

Assume now that $l< k$. From the identity  $H_{E_n}(x^{(n)})=1$ we cannot infer anymore that the sequence $x^{(n)}$ lives in a compact set. Nevertheless, by compactness, we can assume, by taking suitable subsequences, that
\begin{itemize}
\item $p^{(n)}_i(0)$, the $p_i$-th coordinate of $x^{(n)}$, converges to some $\beta_i$, for $1 \le i \le N$.
\item $q^{(n)}_i(0) - q^{(n)}_1(0)$  converges to some $\gamma_i$, for $2\le i\le N$.
\item $E_n^{\unsur{k}-\unsur{l}}q^{(n)}_1(0)$ converges to $\alpha$.
\end{itemize}
Let $\tilde{X}^{(n)}$ be the process of coordinates $(q^{(n)}_i(t) -q^{(n)}_1(0), p^{(n)}_i(t))$. It is starting from $\tilde{x}^{(n)}$ which lies in the hyperplane $\Hrond=\ens{x: q_1=0}$, and converges to some  $\tilde{x}^{(\infty)}\in \Hrond$. The process  $\tilde{X}^{(n)}$ is solution of the system 

\begin{equation}\left\{
 \begin{aligned}
   dq_i(t) &= p_i(t)\, dt& \\
   dp_i(t) &= (-\partial_{q_i}\tilde{H}_{E_n} -\frac{\alen}{2} p_i)\, dt +  E_n^{\unsur{2k} - \frac{3}{4}} \sqrt{T_i} dB_i(t) & (i=1,N)\\
dp_i(t) &= (-\partial_{q_i}\tilde{H}_{E_n})\, dt & \text{\small $(2\le i\le N-1)$}\,.\\
\end{aligned}  \right.
\tag{\text{$\tilde{S}^{E_n}$}}
\end{equation}
 with the Hamiltonian:
$$ \tilde{H}_{E_n}(p,q) = \sum_{1\le i\le N}  \undemi p_i^2 + \unsur{E_n} V(E_n^{\unsur{k}}(q_i -q^{(n)}_1(0))) + \sum_{1\le i\le N-1} \unsur{E_n} U(E_n^{\unsur{k}} (q_{i+1} -q_i))\,.$$

Thanks again to the continuity of solutions of stochastic differential equations with respect to both parameters and starting points we have, for any $\tau>0$:
$$ \lim_{n \to +\infty} P\etc{\sup_{t\le \tau} (\tilde{X}^{(n)}(t) - \tilde{X}^{(\infty)}(t))^2} =0$$
where $\tilde{X}^{(\infty)}$ is the solution  starting from $\tilde{x}^{(\infty)}$ of the system $(S^\infty)$. we can now finish the proof as in the case $k=l$. The only difference is that we have to prove that for a solution of $(S^\infty)$ starting from $x$ such that $H_\infty(x)=1$ and $q_1=0$, we have for any $\tau >0$, 

$$ \int_0^{\tau} (p_1^2(s) + p_N ^2(s))\, ds >0\,.$$
Indeed, if we assume this integral to be $0$, we obtain that on the time interval $\etc{0,\tau}$, the momenta are $0$ and the positions are constant : $q(t)=c$ a vector solution of $\nabla g(c)=0$ if $H_\infty = \sum \undemi p_i^2 + g(c)$, so $c$ is an infimum of the convex function $g$, and since $c_1=0$, this entails $c=0$ which is a contradiction.
\end{proof}

\section{Proof of the main theorem}
\label{sec:proofmain}

We already know, from the preceding sections that the process is irreducible, strong Feller. To conclude the proof of Theorem\ref{sec:introduction:thm}, we shall use Theorem 8.9 of Rey-Bellet\cite{reybellet03}. There only remains to show that there exists $t_0>0$, constants $b_n<+\infty$, $0<\kappa_n< 1$, with $\lim_{n\to +\infty} \kappa_n =0$, and compacts $K_n$ such that
$$ T_{t_0} W(x) \le \kappa_n W(x) + b_n 1_{K_n}(x) $$
If we choose compact sets of the type $K_n=\ens{x: W(x) \le a_n}$ with $a_n \to +\infty$, it is enough to show that 
$$ \lim_{n \to +\infty} \sup_{\ens{x : W(x) > a_n}} \frac{T_{t_0} W(x)}{W(x)} = 0$$
Assume that this is not the case, then we can find $\epsilon >0$ and a sequence $x_n$ with $W(x_n)=e^{\theta H(x_n)} \to +\infty$ and $\frac{T_{t_0} W(x_n)}{W(x_n)}\ge \epsilon$.

Therefore, we shall finish the proof of Theorem~\ref{sec:introduction:thm}, once we have proved 
\begin{lemm}\label{sec:proof-main-theorem:lem}
  Given $0<\theta < \max(T_1,T_N)^{-1}$, there exists a constant $C>0$ and $\alpha,\beta >1$ conjugate exponents, $\unsur{\alpha} + \unsur{\beta}=1$, such that
$$ \frac{T_{t_0} W(x)}{W(x)} \le e^{C \theta (T_1 + T_N) t_0} 
P^x\etc{\exp\etp{ - C \int_0^{t_0} (p_1^2(s) + p_N ^2(s))\, ds}}^{1/\beta}$$
\end{lemm}

Indeed, by combining this Lemma with Lemma~\ref{sec:scaling-1}, we get that
for a subsequence that we call $x_n$ again,
$$ \epsilon \le  e^{C \theta (T_1 + T_N) t_0} \limsup P^{x_n}\etc{\exp\etp{ - C \int_0^{t_0} (p_1^2(s) + p_N ^2(s))\, ds}}^{1/\beta} = 0$$
which is a contradiction.

\begin{proof}[Proof of Lemma~\ref{sec:proof-main-theorem:lem}]
The proof relies on the following clever trick (see e.g. Rey-Bellet and Thomas, proof of Theorem 3.10~\cite{MR1889227}). Since $L$ is the infinitesimal generator of the diffusion $X(t)=(q(t), p(t))$ solution of $(S)$, we have the decomposition
$$ H(X_t) = H(x) + \intot LH(X_s)\, ds + M^H(t)$$
with $M^H(t)= \intot \sigma^t\nabla H(X_s)\, dB_s$ a continuous local martingale of quadratic variation
$$ \crochet{M^H}_t = \intot \norme{\sigma^t\nabla H(X_s)}^2\, ds = 2 \intot \Gamma H(X_s)\, ds$$
where $\Gamma$ denotes the \emph{carré du champ} operator associated to $L$:
$$ \Gamma H(x) = \undemi (L(H^2)(x) - 2 H(x)  L H(x))\,.$$
Only the second order differentials play a role in $\Gamma$ so 
$$ \Gamma H(x) =\undemi \sum_{i=1,N} \frac{T_i}{2} (\partial_{p_i^2}(H^2) - 2 H \partial_{p_i^2} H) = \undemi \sum_{i=1,N} T_i (\partial_{p_i} H)^2 = \undemi \sum_{i=1,N} T_i p_i^2\,.$$
  Hence
  \begin{align}\label{sec:proofmain:1}
LH(x) + \alpha \theta \Gamma H(x) &= \undemi ((T_1 + T_N) - (p_1^2 + p_N^2) + \alpha \theta (T_1 p_1^2 + T_N p_N^2)) \\
&\le \undemi (T_1 + T_N) - \undemi( 1- \alpha \theta \max(T_1,T_N))(p_1^2 + p_N^2). \notag
  \end{align}

Therefore
\begin{align*}
  \frac{T_t W(x)}{W(x)} &=  P^x\etc{ e^{\theta(H(X_t)-H(x))}} = P^x\etc{ e^{\theta M^H_t + \theta \intot LH(X_s)\, ds}}\\
&= P^x\etc{e^{\theta M^H_t - \alpha \frac{\theta^2}{2} \crochet{M^H}_t} e^{\alpha \frac{\theta^2}{2} \crochet{M^H}_t + \theta \intot LH(X_s)\, ds}}\\
&=P^x\etc{U V} \le P^{x}\etc{U^\alpha}^{1/\alpha} P^x\etc{V^{\beta}}^{1/\beta}
\end{align*}
But $U^\alpha$ is just an exponential martingale,
$$ P^{x}\etc{U^\alpha}=P^x\etc{e^{\alpha \theta M^H_t - \undemi \alpha^2 \theta^2 \crochet{M^H}_t}} = 1$$
Hence, injecting the inequality~\eqref{sec:proofmain:1}, we obtain
\begin{align*}
\frac{T_t W(x)}{W(x)} &\le e^{\undemi t \beta \theta (T_1 + T_N)} \times \ldots \\
&\times P^{x}\etc{e^{-\frac{\beta\theta}{2} (1 - \alpha \theta \max(T_1, T_N)) \intot (p_1^2(s) + p_2^2(s))\, ds}}^{1/\beta}
\end{align*}

If we choose $\alpha$ close enough to $1$ so that $\alpha \theta \max(T_1, T_N)<1$, we obtain the desired inequality.
\end{proof}

\section{The moment exchanging model}

The Hamiltonian is harmonic with no pinning (see~\cite{MR2185330})
$$ H = \undemi \sum_{1\le i\le N}{p_i^2} + \undemi \sum_{1\le i \le N-1} (q_{i+1} -q_i)^2\,.$$
Considering the interdistances $r_i=q_{i+1} -q_i$, we not only have two heat baths at both ends of the chain, but also a momentum exchanging noise. The stochastic dynamics are described by the stochastic differential system:

\begin{equation}\left\{
 \begin{aligned}
   dr_i(t) &= (p_{i+1} -p_i)\, dt\\
   dp_1(t) &= (r_1 -\frac{1+\gamma}{2} p_1)\, dt -\sqrt{\gamma} p_2 dB_{1,2}(t) + \sqrt{T_1} dB_{0,1}(t) & \\
   dp_N(t) &= (-r_{N-1} - \frac{1+\gamma}{2}p_N)\, dt +  \sqrt{\gamma} p_{N-1} dB_{N-1,N}(t) + \sqrt{T_N} dB_{N,N+1}(t) & \\
dp_i(t) &= (r_i -r_{i-1}-\gamma p_i)\, dt + \sqrt{\gamma} (p_{i-1} dB_{i-1,i}(t) - p_{i+1} dB_{i,i+1}(t))& \text{\small $(2\le i\le N-1)$}
\end{aligned}  \right.
\tag{\text{$\Sigma$}}
\end{equation}
where $(B_{j,j+1})_{0\le i\le N}$ are independent Brownian motions. The infinitesimal generator is now
  \begin{gather*}
L = \Arond + \Lrond_R + \frac{\gamma}{2} \Srond \\
\Srond = \sum_{1\le i\le N-1} X_{x,x+1}^2\,, \quad X_{i,i+1} = p_{i+1} \partial_{p_i} - p_i \partial_{p_{i+1}}\,.    
  \end{gather*}

We have : $ X_{i,i+1}(p_i^2) = 2 p_i p_{i+1} = -  X_{i,i+1}(p_{i+1}^2)$. Therefore, $X_{i,i+1} H = 0$, $\Srond(H)=0$ and 
$$ \Srond(W) = \sum_{i}W (\theta^2 (X_{i,i+1}H)^2 + \theta X^2_{i,i+1}(H)) = 0\, .$$ This implies that we can proceed through section 2 with no change. In particular:
$LW \le \undemi \theta (T_1 + T_N) W$ if $0 <\theta \max(T_1,T_N) < 1$.

\medskip

 Since we have more squared vector fields, we also have H\"ormander's hypoellipticity condition. We also dispose of more control, so the semigroup is irreducible (section 4).

\medskip

To use the scaling technique of section 5, we need to look closely at the semimartingale decomposition of $H(X_t)$: 
$$ H(X_t) = H(x) + M^H_t + \intot LH(X_s)\, ds$$
with
\begin{align*}
  dM^H_t &= p_1 (-\sqrt{\gamma} p_2 dB_{1,2}(t) + \sqrt{T_1} dB_{0,1}(t)) \\
&+ \sum_{2\le i \le N-1} \sqrt \gamma p_i ( p_{i-1} dB_{i-1,i} - p_{i+1} dB_{i,i+1}) \\
& + p_N (\sqrt{\gamma} p_{N-1} dB_{N-1,N} + \sqrt{T_N} dB_{N,N+1} \\
&= p_1\sqrt{T_1} dB_{0,1}(t) + p_N dB_{N,N+1}(t)\,.
\end{align*}
Therefore, $d\crochet{M^H}_t= (p_1^2 T_1 + p_N^2 T_N)dt$ and the proofs of section 6 can proceed (almost) unchanged. The only thing that changes is that the limiting system is no more deterministic but stochastic : only the heat reservoir noises disappear

\begin{equation}\label{bolla:sysinfini}
\left\{
 \begin{aligned}
   dr_i(t) &= (p_{i+1} -p_i)\, dt\\
   dp_1(t) &= (r_1 -\frac{1+\gamma}{2} p_1)\, dt -\sqrt{\gamma} p_2 dB_{1,2}(t)  & \\
   dp_N(t) &= (-r_{N-1} - \frac{1+\gamma}{2}p_N)\, dt +  \sqrt{\gamma} p_{N-1} dB_{N-1,N}(t) & \\
dp_i(t) &= (r_i -r_{i-1}-\gamma p_i)\, dt + \sqrt{\gamma} (p_{i-1} dB_{i-1,i}(t) - p_{i+1} dB_{i,i+1}(t))& \text{\small $(2\le i\le N-1)$}
\end{aligned}  \right.
\tag{\text{$\Sigma^\infty$}}
\end{equation}

We only need to replace Lemma~\ref{sec:scaling:lem} with the following non degeneracy result:
\begin{lemm}\label{sec:bolla:lem}
  Let $(q(t),p(t))$ be a solution of $(\Sigma^\infty)$, starting from $x$ such that $H_\infty(x)=1$. Then, for any $\tau >0$, almost surely,
$$ \int_0^\tau (p_1^2(s) + p_N ^2(s))\, ds >0\,.$$
\end{lemm}
\begin{proof}
  Assume on the contrary that 
$$ \int_0^\tau (p_1^2(s) + p_N ^2(s))\, ds =0\,\quad a.s.$$
Then on $I=\etc{0,\tau}$, we have $p_1=0$ a.s., and thus
$$ \intot r_1(s)\, ds = \intot \sqrt{\gamma} p_2(s) dB_{1,2}(s)\,,\quad (\forall t\le \tau)$$
This states that a finite variation process coincides with a continuous martingale. Hence,  they both vanish : on $I$, $r_1=p_2=0$ a.s. It should be now clear how to proceed by induction and to obtain $r_i=0=p_i$, which contradicts the fact that $H(p,q)=1$ (here $H_\infty=H$).
\end{proof}

\section{The Lefevere-Schenkel chain}
Lefevere and Schenkel~\cite{1742-5468-2006-02-L02001} consider a \emph{periodic} lattice Hamiltonian 
$$ H = \undemi \sum_{i=1}^N p_i^2 + \omega^2 \mu^2 q_i^2 + \omega^2 (q_i - q_{i-1})^2 $$
 (we have $p_{N+k}=p_k, q_{N+k}=q_k$ and for sake of notations $N$ is a multiple of $4$). They introduce the Fourier coordinates
$$ Q_k = \unsur{\sqrt{N}} \sum_{i=1}^N e^{i \frac{2\pi k}{N} j} q_j\,, \quad 
P_k = \unsur{\sqrt{N}} \sum_{i=1}^N e^{i \frac{2\pi k}{N} j} p_j\,,$$
which satisfy $P_{-k} = P_k^*$ (the complex conjugate) and $Q_{-k}=Q_k^*$. This implies in particular that $Q_0,P_0,Q_{N/2},P_{N/2}$ are real valued since for example $P_{N/2} = P_{N/2 -N} = P_{-N/2}=P_{N/2}^*$. In Fourier coordinates the Hamiltonian reads:
$$ H = \undemi \sum_{-\frac{N}{2} + 1 \le k \le \frac{N}{2}}\valabs{P_k}^2 + \omega_k^2 \valabs{Q_k}^2$$
with $\valabs{z}^2=zz^*$ the complex square modulus and $\omega_k^2=\omega^2(\mu^2 + 4 \sin^2(\frac{k\pi}{N}))$.

The momenta $P_k$ are coupled to heat baths at a fixed constant temperature $T>0$, whereas the positions $Q_k$ are coupled to heat baths at temperatures $D_k^2$ where to follow notations of Lefevere and Schenkel~\cite{1742-5468-2006-02-L02001} we have $D_k = \unsur{T} \tau^2 \alpha^2(\frac{k \pi}{N})$. In short they are solutions to the system of stochastic differential equations
\begin{equation}\left\{
 \begin{aligned}
dP_k(t) &= -(\omega_k^2 Q_k +\undemi P_k)\, dt + \sqrt{T} dZ_k(t) \\
dQ_k(t) &= P_k  dt - i D_k dZ_k(t)
\end{aligned}  \right.
\tag{\text{$LS$}}
\end{equation}

where $(Z_k(t))_{t\ge 0}$ are independent standard complex Brownian motions if $k\neq 0,N/2$ and ordinary real valued Brownian motions if $k=0,N/2$. The boundary conditions imply that $D_{k}=D_{-k}=D_{k+N}$ (and thus $D_0=D_{N/2}=0$) and that $Z_{-k}(t)= Z_k^*(t)$.

We shall work with real valued processes : $P_k= R_k + i S_k$, $Q_k = V_k + i W_k$, $Z_k = X_k + i Y_k$ for $k \neq 0, N/2$. The stochastic differential system is for $k\neq 0, N/2$:

\begin{equation}\left\{
 \begin{aligned}
dR_k(t) &= -(\omega_k^2 V_k +\undemi R_k)\, dt + \sqrt{T} dX_k(t) \\
dS_k(t) &= -(\omega_k^2 W_k +\undemi S_k)\, dt + \sqrt{T} dY_k(t) \\
d V_k(t) &= R_k dt + D_k dY_k(t)\\
dW_k(t) &= S_k dt - D_k dX_k(t)
\end{aligned}  \right.
\tag{\text{$LS_4$}}
\end{equation}

and for $k=0,N/2$ ( $D_k=0$ and $P_k,Q_k,Z_k$ are real valued):

\begin{equation}\left\{
 \begin{aligned}
dR_k(t) &= -(\omega_k^2 V_k +\undemi R_k)\, dt + \sqrt{T} dX_k(t)\,, \\
d V_k(t) &= R_k dt \,.\\
\end{aligned}  \right.
\tag{\text{$LS_2$}}
\end{equation}

Therefore the infinitesimal generator can be written $\Lrond = \Arond + \Lrond_R$ with $\Arond$ the Hamiltonian generator and
\begin{align*}
  \Lrond_R &= \sum_{p=R_k,S_k} \undemi(-p \partial_p + T \partial^2_{p^2}) \\
&+ \sum_k \undemi D_k^2 (\partial^2_{V_k^2} +\partial^2_{W_k^2} ) \\
&+ \sum_k \sqrt{T} D_k (\partial^2_{S_k V_k} -\partial^2_{R_k W_k})\,.
\end{align*}
H\"ormander's hypoellipticity condition is very easy to check and we have
$$ \Lrond H = \Arond H + \Lrond_R H = 0 + NT+ \sum_k D_k^2\omega_k^2 - \undemi \sum_k \valabs{P_k}^2\le C\,,$$
with $C$ a constant, which we can choose  $C\ge 1$, and with the function $W=C+H$ that satisfies $LW \le W$, we see that there is no explosion.

The control problem is easy to solve since it splits into two elementary control problems, one in $\R^4$ and one in $\R^2$. Let us examine the control problem in $\R^4$: $\frac{dx}{dt} = A x + B u(t)$, with
$$ A = 
\begin{pmatrix} -\undemi & 0 & -\omega_k^2 & 0 \\
0 & -\undemi & 0 & -\omega_k^2 \\
1 & 0 & 0 & 0 \\
0 & 1 & 0 & 0\\
\end{pmatrix}
\quad B= 
\begin{pmatrix}
\sqrt{T} & 0 \\
0 & \sqrt{T} \\
0 & D_k \\
-D_k & 0  
\end{pmatrix}
$$
We let $\Srond$ be the smallest vector space containing the range of $B$ and stable by $A$. $\Srond$ contains the vector $f_1=( \sqrt{T}, 0, 0, -D_k)^T$ and $Af_1$ so it contains $A f_1 + \undemi f_1 = (-\omega_k^2 D_k, 0, 0, \sqrt{T})^T$ and thus $\Srond$ contains the coordinate vectors $(1,0,0,0)^T$ and $(0,0,0,1)^T$. Working similarly with $f_2=(0,\sqrt{T},D_k,0)^T$ we conclude that $\Srond$ is $\R^4$ and by  Kalman's criterion the accessibility sets are both $\R^4$.

To use the scaling technique of section~\ref{sec:scaling}, we need to take a close look at the semimartingale decomposition of $H(\xi_t)$ with $\xi_t=(P_k(t),Q_k(t),\frac{N}{2}+1 \le k\le \frac{N}{2})$.
$$ dH(\xi_t) = dM^H_t + \Lrond H(\xi_t)\, dt $$
with $M^H_t$ a continuous local martingale:
$$ dM^H_t = \sum_k (\sqrt{T} -\omega_k^2 D_k W_k) dX_k(t) + (\sqrt{T} S_k + \omega_k^2 D_k V_k) \, dY_k(t)\,,$$
with continuous quadratic variation $d\crochet{M^H}_t = 2 \Gamma H(\xi_t)\, dt$. Hence:
\begin{align*}
  \Gamma H(\xi) &= \undemi \sum_k  (\sqrt{T} -\omega_k^2 D_k W_k)^2 + (\sqrt{T} S_k + \omega_k^2 D_k V_k)^2 \\
& \le \sum_k T \valabs{P_k}^2 + \omega_k^4 D_k^2 \valabs{Q_k}^2
\end{align*}

Proceeding as in section~\ref{sec:proofmain}, we see that
$$ LH(\xi) + \theta \alpha \Gamma H(\xi) \le C + \sum_k (-\undemi + \alpha \theta T) \valabs{P_k}^2 + \alpha \theta \omega_k^4 D_k^2 \valabs{Q_k}^2$$
Hence, if we fix $\alpha>1$, there exists constants $\theta_0,C_1,C_2,C_3>0$ such that for $0< \theta < \theta_0$,
$$ LH(\xi) + \theta \alpha \Gamma H(\xi) \le C_1 - \undemi C_2\sum_k \etp{ \valabs{P_k}^2 - \theta C_3 \valabs{Q_k}^2}\,.$$
This yields an analogue of Lemma~\ref{sec:proof-main-theorem:lem} and the proof can proceed as in section~\ref{sec:proofmain}. More precisely, we have the upper bound for $0< \theta < \theta_0$ , $W=e^{\theta H}$ 

\begin{equation}
  \label{eq:lefevere:majsg}
  \frac{T_t W(x)}{W(x)} \le e^{C_1 t} P^x\etc{\exp -\frac{C_2}{2} \intot \sum_k \etp{ \valabs{P_k}^2 - \theta C_3 \valabs{Q_k}^2}(s)\, ds}^{1/\beta}
\end{equation}

 There is only one thing left to establish, the non degeneracy of the limiting system : the noise disappear, so that $(R_k,V_k)$ and $(S_k,W_k)$ are solutions of the same deterministic system, with $\omega^2=\omega_k^2$,

\begin{equation}\left\{
 \begin{aligned}
dr &= -(\omega^2 v + \undemi r)\, dt  \\
dv &= r dt\\
\end{aligned}  \right.
\tag{\text{$L_\infty$}}
\end{equation}

\begin{lemm}
  Let $x=(r,v)$ be a solution of $L_\infty$ with starting point $x_0\neq(0,0)$. Then
$$ \liminf_{\tau \to +\infty} \frac{\int_0^\tau r^2(s)\, ds}{\int_0^{\tau} v^2(s)\, ds} >0\,.$$
\end{lemm}
\begin{proof}
  This is a simple exercise in ordinary differential equations. The solution of the system is $x(t)= e^{tA}x_0$. The eigenvalues of $A$ are
$$ \lambda_{\pm} = 
\begin{cases}
  -\unsur{4} \pm \sqrt{\unsur{16} - \omega^2} &\text{if $\omega^2 \le \unsur{16}$;} \\
 -\unsur{4} \pm i\sqrt{ \omega^2-\unsur{16}}&\text{otherwise.}
\end{cases}
$$
Hence, $v=A_+ e^{\lambda_+ t} + A_- e^{\lambda_-} t$, $r(t) = \frac{dv}{dt}$ and thus, in the first case, $\omega^2 \le \unsur{16}$ :
$$ \lim_{\tau \to +\infty} \frac{\int_0^\tau r^2(s)\, ds}{\int_0^{\tau} v^2(s)\, ds} =
\begin{cases}
  \lambda_+^2 &\text{if $A_+\neq 0$;} \\
\lambda_-^2 &\text{otherwise.}
\end{cases}
$$
We leave the second case  $\omega^2 > \unsur{16}$ to the interested reader.
\end{proof}

\bibliographystyle{amsplain}
\bibliography{new,poly}

\end{document}